\documentclass[letter]{amsart}
\usepackage{amsmath,amsfonts,amsthm,amssymb}
\usepackage{hyperref} 
\usepackage{graphicx}
\usepackage{epstopdf}

\theoremstyle{plain}
\newtheorem{theorem}{Theorem}[section]

\newtheorem{lemma}[theorem]{Lemma}
\newtheorem*{de-lemma}{Lemma}
\newtheorem{proposition}[theorem]{Proposition}

\theoremstyle{remark}

\theoremstyle{definition}

\newcommand{\R}{\mathbb{R}}

\begin{document}

\title[Blowing up solutions of semilinear PDE with convex potentials]{Blowing up solutions of semilinear P.D.E. with convex potentials}

\author{Panayotis Smyrnelis} \address[P.~ Smyrnelis]{Institute of Mathematics,
Polish Academy of Sciences, 00-956 Warsaw, Poland}
\email[P. ~Smyrnelis]{psmyrnelis@impan.pl}

\date{}

\maketitle
\begin{abstract}
We consider convex potentials $W:\R\to [0,\infty)$ vanishing at $0$ and growing sufficiently fast at $\pm\infty$. Given any open set $\Omega\subset\R^n$ with Lipschitz and compact boundary, we prove the existence and uniqueness of a solution of $\Delta u= W'(u)$ in $\Omega$, such that $u=+\infty$ or $u=-\infty$ on $\partial \Omega$. Moreover, if $\partial \Omega$ is the union of two disjoint compact subsets $A^+$ and $A^-$, there also exists a unique solution satisfying $u=+\infty$ on $A^+$ and $u=-\infty$ on $A^-$. 

\end{abstract}

\section{Introduction and statement of the main theorem}

In this paper we study the solutions of the P.D.E.
\begin{equation}\label{pde}
\Delta u=W'(u), \ u:\R^n \supset \Omega \to \R, 
\end{equation}
where $\Omega\subset\R^n$ is an open set (possibly unbounded) with Lipschitz boundary, and $W:\R\to \R$ is a potential satisfying the following hypotheses:
\begin{subequations}\label{hhhh}
\begin{equation}\label{hh1}
\text{$W\in C^2(\R,\R)$ is nonnegative, convex, and vanishes only at $0$,}
\end{equation}
\begin{equation}\label{growW}
\int_{|u|>1}\frac{du}{\sqrt{W(u)}}<\infty.
\end{equation}
\end{subequations}
For instance we can take $W(u)=|u|^\alpha$, with $\alpha>2$, or $W(u)=\cosh u-1$.
The corresponding O.D.E. 
\begin{equation}\label{ode}
u''=W'(u), \ u:\R \supset I  \to \R, 
\end{equation}
is Newton's one dimensional equation for a unit mass and \emph{potential energy} $-W$, 
where the variable $x \in \R$ stands for time. We recall that the \emph{Hamiltonian} 
\begin{equation}\label{ham}
H:=\frac{1}{2}|u'|^2-W(u),
\end{equation}
or total mechanical energy, is constant along solutions. Assuming that $W$ satisfies \eqref{hhhh}, the maximal solutions of \eqref{ode} such that $H\neq 0$ are defined on bounded intervals (cf. Proposition \ref{prop1}). More precisely, for every bounded interval $I=(a,b)$, there exists a unique solution of \eqref{ode} diverging to $+\infty$ (resp. $-\infty$) at the endpoints $a$ and $b$. In addition, there also exists a unique solution such that $u=-\infty$ at $a$, and $u=+\infty$ at $b$ (resp. $u=+\infty$ at $a$, and $u=-\infty$ at $b$). On the other hand, $0$ is the only solution defined on all $\R$, while the remaining solutions with $H=0$ are defined on half-bounded intervals $(a,+\infty)$ or $(-\infty,b)$: they converge to $0$ at $+\infty$ (resp. $-\infty$), and blow up at $a$ (resp. $b$). In the case where instead of \eqref{growW}, we assume that $\int_{u>1}\frac{du}{\sqrt{W(u)}}=\infty$, and $\int_{u<-1}\frac{du}{\sqrt{W(u)}}=\infty$, then every maximal solution of \eqref{ode} is defined on $\R$.

The scope of this paper is to construct in a general setup the analog of the aforementioned orbits for P.D.E. \eqref{pde}, that is, solutions which are infinite on $\partial \Omega$: cf. (i) and (ii) in Theorem \ref{th1} below. Theorem \ref{th1} (i) has first been proved in \cite{nirenberg} for the special choice of $W(u)=|u|^{\frac{2n}{n-2}}$, in order to investigate certain conformally invariant P.D.E. associated with Riemannian metrics. 
For more applications in geometric analysis, in particular in the theory of mean curvature surfaces and harmonic maps, we refer to \cite{han}, \cite{wan}, \cite{guan}, \cite{aviles}, and \cite{anestis}.
In contrast with \cite{nirenberg}, where $\partial \Omega$ is assumed to be the union of smooth manifolds of sufficiently small codimension, we consider open sets $\Omega$ with Lipschitz boundary. As far as we know, the solutions in Theorem  \ref{th1} (ii) have not appeared in the litterature.

\smallskip

\begin{theorem}\label{th1}
Given any open set $\Omega\subset\R^n$ with Lipschitz boundary, and assuming \eqref{hhhh}:
\begin{itemize}
\item[(i)] There exists a positive solution $u\in C^2(\Omega)$ of \eqref{pde} such that
\begin{equation}\label{bc1}
\forall x_0\in\partial\Omega: \ \lim_{\Omega\ni x \to x_0} u(x)=+\infty \text{ (resp. $-\infty$)}, 
\end{equation}

\item[(ii)] If $\partial \Omega$ is partitioned into two disjoint closed subsets $A^-$ and $A^+$, there also
exists a solution $u\in C^2(\Omega)$ of \eqref{pde} such that
\begin{equation}\label{bc3}
\forall x_0\in A^+ \text{ (resp. $A^-$)}: \ \lim_{\Omega\ni x\to x_0} u(x)=+\infty \text{ (resp. $-\infty$)}.\footnote{In view of \eqref{bc3}, this solution `connecting' $+\infty$ to $-\infty$, may be compared to the \emph{heteroclinic} orbits that connect in a similar way two distinct equilibrium points (cf. for instance \cite{indiana}).}
\end{equation}
\item[(iii)] In addition, denoting by $d$ the Euclidean distance, the solutions in (i) and (ii) above satisfy
\begin{equation}\label{bc2}
\lim_{x\in\Omega,\ d(x,\partial \Omega)\to\infty} u(x)=0, \text{ provided that the function $\Omega\ni x\mapsto d(x,\partial\Omega)$ is unbounded.} 
\end{equation}

\item[(iv)] Assuming moreover that $\partial \Omega$ is compact, and $W'$ is convex on $(0,\infty)$ and concave on $(-\infty,0)$, the solution of \eqref{pde} satisfying \eqref{bc1} (resp. \eqref{bc3}) is unique.
\end{itemize}
\end{theorem}

\smallskip

\section{Structure of solutions of O.D.E. \eqref{ode}}
For the convenience of the reader we recall in this section the structure of solutions of O.D.E. \eqref{ode}.
In the phase plane these orbits are described by equation \eqref{ham}, expressing the conservation of the total mechanical energy. 
Depending on the sign of $H$, we present in Proposition \ref{prop1} below, the four different kinds of orbits obtained. 
The solutions with negative Hamiltonian $H=-W(\lambda)<0$, for some $\lambda\neq 0$, correspond to case (a), while the ones with positive $H=h>0$, to case (b). Finally, the solutions whose Hamiltonian vanishes are examined in cases (c) and (d). 
\smallskip
\begin{proposition}\label{prop1}  
Assuming \eqref{hhhh}, any solution of \eqref{ode} coincides up to translation and change of $x$ by $-x$, with one of the solutions $\alpha$, $\beta$, $\gamma$, $\tilde \gamma$, or $0$ described below:
\begin{itemize}
\item[(a)] For every $\lambda>0$, the maximal solution $\alpha$ of \eqref{ode} such that \begin{equation}\label{alpha}
\alpha(0)=\lambda \text{ and } \alpha'(0)=0, 
\end{equation}
is defined on the interval $(-l,l)$, with 
\begin{equation*}
l(\lambda)=\int_0^\infty\frac{ds}{\sqrt{2(W(s+\lambda)-W(\lambda))}}\in (0,\infty).
\end{equation*}
$\alpha$ is also even, and $\alpha:[0,l)\to[\lambda,+\infty)$ is strictly increasing and onto. Moreover, the function $(0,\infty)\ni\lambda\mapsto l(\lambda)\in (0,\infty)$ is continuous, strictly decreasing, and onto.

Similarly when $\lambda<0$, the maximal solution of \eqref{ode} satisfying \eqref{alpha} is even and defined on the interval $(-l,l)$ with
\begin{equation*}
l(\lambda)=\int_{-\infty}^0\frac{ds}{\sqrt{2(W(s+\lambda)-W(\lambda))}}\in (0,\infty).
\end{equation*}
In addition, $\alpha:[0,l)\to (-\infty,\lambda]$ is strictly decreasing and onto, while $(-\infty,0)\ni l\mapsto l(\lambda)\in (0,\infty)$, is continuous, strictly increasing and onto.
\item[(b)] For every $h>0$, the maximal solution $\beta$ of \eqref{ode} such that $\beta(0)=0$, and $\beta'(0)=\sqrt{2h}$, is defined on the interval $(-l^-,l^+)$, with 
$$l^+(h)=\int_0^\infty\frac{ds}{\sqrt{2(W(s)+h)}}\in (0,\infty),$$
and $$l^-(h)=\int^0_{-\infty}\frac{ds}{\sqrt{2(W(s)+h))}}\in (0,\infty).$$
$\beta:(-l^-,l^+)\to (-\infty,+\infty) $ is also strictly increasing, and onto. Moreover, the functions $(0,\infty)\ni h \mapsto l^\pm(h)\in (0,\infty)$ are continuous, strictly decreasing, and onto. 
\item[(c)] There exists a unique solution $\gamma$ (resp. $\tilde \gamma$) of \eqref{ode} defined on $(0,\infty)$ and such that $\lim_0\gamma=\infty$ (resp. $\lim_0\tilde \gamma=-\infty$). Moreover, $\gamma:(0,\infty)\to (0,\infty)$ (resp. $\tilde\gamma:(0,\infty)\to (-\infty,0)$) is strictly decreasing (resp. increasing), and onto.
\item[(d)] $u\equiv 0$ is the only solution of \eqref{ode} defined on all $\R$.
\end{itemize}
\end{proposition}
\emph{Proof. }
(a) A maximal solution $u:\R\supset (a,b)\to \R$ of \eqref{ode} such that $H<0$ satisfies $W(u)\geq -H>0$ in view of \eqref{ham}. Thus, if $\tilde \lambda<0$ and $\lambda>0$ are the two roots of the equation $W(t)=-H$, we will have either $u(I)\subset [\lambda,\infty)$ or $u(I)\subset (-\infty,\lambda']$. In what follows, we shall examine only the former case, since the latter is similar. Actually, since the maximal solution $u$ is convex on $(a,b)$, we have $u((a,b))=[\lambda,\infty)$, and $u(x)\to \infty$, as $x\to a$ (resp. $x\to b$). Thus, $u(x_0)=\lambda$, and $u'(x_0)=0$, hold for some $x_0\in (a,b)$. This implies that up to a translation $u$ coincides with $\alpha$. In addition, since $x\mapsto \alpha(-x)$ is also a solution of \eqref{ode} satisfying \eqref{alpha}, we can see that $\alpha $ is even. Let us now compute the length of the interval $(-l,l)$ where $\alpha$ is defined. The derivative $\alpha'$ is positive on $(0,l)$, hence $\alpha:[0,l)\to[\lambda,\infty)$ is strictly increasing and onto. As a consequence of \eqref{ham}, $\alpha'(x)=\sqrt{2(W(\alpha(x))-W(\lambda))}$ for $x> 0$. Setting $\phi=\alpha^{-1}:(\lambda,\infty)\to(0,l)$, we have $\phi'(\alpha)=\frac{1}{\sqrt{2(W(\alpha)-W(\lambda))}}$, and in view of \eqref{growW}:
$$l=\int_{\lambda}^{\infty}\phi'(\alpha)d\alpha=
\int_0^\infty\frac{ds}{\sqrt{2(W(s+\lambda)-W(\lambda))}}< \infty.$$
Moreover, since $W$ is convex, one can see that $\lambda_1<\lambda_2$ implies that 
\begin{equation}
\label{notall}
\frac{1}{\sqrt{2(W(s+\lambda_2)-W(\lambda))}}\leq \frac{1}{\sqrt{2(W(s+\lambda_1)-W(\lambda))}}, \ \forall s\in [0,\infty),
\end{equation}
and due to \eqref{growW} the equality in \eqref{notall} cannot hold for every $s \geq 0$. Thus, $l(\lambda_2)<l(\lambda_1)$. From \eqref{notall} it also follows by dominated convergence that $(0,\infty)\ni\lambda\mapsto l(\lambda)$ is continuous, and that $\lim_{\lambda\to\infty}l(\lambda)=0$, since 
 $\lim_{\lambda \to \infty}\frac{1}{\sqrt{2(W(s+\lambda)-W(\lambda))}}=0$, $\forall s>0$, in view of \eqref{growW}. Finally, we obtain by monotone convergence  $\lim_{\lambda\to 0}l(\lambda)=\infty$, since $$\lim_{\lambda \to 0}\int_0^\infty\frac{ds}{\sqrt{2(W(s+\lambda)-W(\lambda))}}=\int_0^\infty\frac{ds}{\sqrt{2W(s)}}=\infty.$$

(b) Similarly, one can show that a maximal solution $u:\R\supset (a,b)\to\R$ of O.D.E. \eqref{ode}, with $H>0$, is strictly monotone and such that $u((a,b))=\R$. Up to translation, and change of $x$ by $-x$, it coincides with the solution $\beta$. The length of the interval $(-l^-,l^+)$ where $\beta$ is defined, is determined as in (a), and the properties of $l^\pm$ are established in a similar way.

(c) and (d) Finally, a maximal solution of \eqref{ode} such that $H=0$, and $u(x_0)=0$ for some $x_0\in \R$, is identically $0$ in view of \eqref{ham} and the uniqueness result for O.D.E. One can easily prove that the remaining maximal solutions such that $H=0$, coincide up to translation and change of $x$ by $-x$, either with the solution $\gamma$ or $\tilde \gamma$.
\qed

\section{Construction of supersolutions and subsolutions}
We shall utilize two basic results for P.D.E. \eqref{pde}: the existence of a unique solution of the Dirichlet problem in a bounded domain with bounded boundary condition, and the maximum principle.
\smallskip
\begin{lemma}\label{lem01}  
Let $W$ be a potential satisfying \eqref{hh1}, and let $\omega\subset  \R^n$ be a bounded open set with Lipschitz boundary.
Then, given $\phi\in W^{1,2}(\omega)\cap L^\infty(\omega)$, the Dirichlet problem
\begin{equation}\label{Diri}
\Delta u =W'(u) \text{ in $\omega$, with $u=\phi$ on $\partial \omega$ (in the sense of the trace),}
\end{equation}
has a unique solution in $W^{1,2}(\omega)\cap L^\infty(\omega)$.
\end{lemma}
\emph{Proof. } The solution of \eqref{Diri} is the global minimizer of the energy functional $E(u)=\int_\omega\big( \frac{1}{2}|\nabla u|^2+W(u)\big)$ in the closed affine subspace $\mathcal A:=\{u\in W^{1,2}(\omega): u-\phi\in W^{1,2}_0(\omega)\}$. The uniqueness of the solution follows by the strict convexity of $E$.
\qed
\bigskip
\begin{lemma}\label{lem02}  
Let $W$ be a potential satisfying \eqref{hh1}, and let $\omega\subset  \R^n$ be a bounded open set with Lipschitz boundary.
Then, given $u_1,u_2\in W^{1,2}(\omega)\cap L^\infty(\omega)$, such that
\begin{equation*}
\Delta u_1\geq W'(u_1) \text{ and }\Delta u_2\leq W'(u_2) \text{ hold weakly in $\omega$,}
\end{equation*}
the condition $u_1\leq u_2$ on $\partial \omega$ (in the sense of the trace) implies that $u_1\leq u_2$ a.e. in $\omega$.
\end{lemma}
\emph{Proof. } Setting $v=u_1-u_2$, we have by the convexity of $W$:
\begin{equation*}
\Delta v(x)\geq W'(u_1(x))-W'(u_2(x))=c(x) v(x),
\end{equation*}
with $c(x):=\int_0^1W''(u_1(x)+tv(x))dt\geq 0
$, and $c\in L^\infty(\omega)$. Thus, the maximum principle can be applied to $v$ to deduce that $v\leq 0$.
\qed
\bigskip

From Proposition \ref{prop1} we easily obtain a family of supersolutions and subsolutions of \eqref{pde}: \smallskip
\begin{lemma}\label{lem1}  
Let $\tilde \alpha$ be the orbit (cf. \eqref{alpha}) provided by Proposition \ref{prop1} (a) for the potential $\tilde W=\frac{1}{n}W$. Then,
$\zeta(x)=\tilde \alpha_\lambda (|x|)$ is a supersolution (resp. subsolution) of \eqref{pde} when $\lambda>0$ (resp. $\lambda<0$). As a consequence, $0$ is the only solution of \eqref{pde} defined in $\R^n$.
\end{lemma}

\emph{Proof. } We first show that $\zeta$ is a supersolution when $\lambda>0$.
Indeed, by the mean value theorem, and the convexity of $W$, we have:
\begin{align*}
\Delta \zeta(x)&=\tilde\alpha'' (|x|)+\frac{n-1}{|x|}\tilde\alpha'(|x|)= \frac{1}{n}W'(\tilde \alpha(|x|)+(n-1)\tilde\alpha''(\xi), \text{ for some $\xi\in (0,|x|)$,}\\
&= \frac{1}{n}W'(\tilde\alpha(|x|)+\frac{n-1}{n}W'(\tilde \alpha(\xi))\leq W'(\zeta(x)), \text{ provided that $\lambda>0$}.
\end{align*}
In the same way, we establish that $\zeta$ is a subsolution when $\lambda<0$.
Finally, if $u\in C^2(\R^n)$ is a solution of \eqref{pde}, we can see that $u$ is bounded above by the functions $\tilde \alpha_\lambda(|x-x_0|)$, for every $x_0\in \R^n$, and $\lambda >0$. Indeed, it is clear that $u(x)\leq \zeta(x-x_0)$ holds for $|x-x_0|\geq l(\lambda)-\epsilon$, with $\epsilon>0$ small enough. Thus, in view of Lemma \ref{lem02}, $u(x)\leq \zeta(x-x_0)$ also holds for $|x-x_0|< l(\lambda)-\epsilon$. Letting $\lambda\to 0$, we deduce that $u\leq 0$. The proof that $u\geq 0$, is similar.  
\qed

\bigskip
Next, in order to establish the boundary conditions \eqref{bc1} and \eqref{bc3}, we construct `barrier' functions defined in annuli. Property (ii) in Lemma \ref{lem2} is essential to address the issue of Lipschitz boundaries. In what follows we denote by $B(x_0,r)$ the ball of radius $r$ centered at $x_0$.
\smallskip
\begin{lemma}\label{lem2}
Given $L>1$, there exists for every $\epsilon>0$, a radial solution $x\mapsto\xi_\epsilon(x)=\tilde\xi_\epsilon(|x|)$ of \eqref{pde} defined in the annulus $\overline{B(0,L\epsilon)}\setminus \overline{B(0,\epsilon)}$, and such that 
\begin{itemize}
\item[(i)] $\tilde \xi_\epsilon \in C^2((\epsilon,L\epsilon])$ is strictly decreasing, with $\lim_{r\to \epsilon}\tilde\xi_\epsilon(r)=\infty$, and $\tilde\xi_\epsilon(L\epsilon)=0$,
\item[(ii)] $\lim_{\epsilon\to 0}\tilde\xi_\epsilon (r \epsilon)=\infty$, holds for every $r \in (1,L)$ fixed.
\end{itemize}  
\end{lemma}

\emph{Proof. } It is convenient to work in a fixed annulus $\omega=B(0,L)\setminus \overline{ B(0,1)}$, by setting $\psi_\epsilon(x)=\xi_\epsilon(\epsilon x)$, and to seek instead of $\xi_\epsilon$, a radial solution of
\begin{equation}\label{Dir2}
\Delta \psi_\epsilon =\epsilon^2W'(\psi_\epsilon) \text{ in $\omega$, such that $\lim_{|x|\to1}\psi_\epsilon(x)=\infty$, and  $\psi_\epsilon(x)=0$ for $|x|=L$.}
\end{equation}
We first consider the solution $\psi_{\epsilon,n}\in C^2(\overline \omega)$ of the Dirichlet problem
\begin{equation}\label{Dir3}
\Delta \psi_{\epsilon,n} =\epsilon^2W'(\psi_{\epsilon,n}) \text{ in $\omega$, with $\psi_{\epsilon,n}(x)=n$ for $|x|=1$, and  $\psi_{\epsilon,n}(x)=0$ for $|x|=L$.}
\end{equation}
It is easy to see that $\psi_{\epsilon,n}$ is radial i.e. $\psi_{\epsilon,n}(x)=\tilde\psi_{\epsilon,n}(|x|)$, with $\tilde \psi_{\epsilon,n}$ strictly decreasing and convex on $[1,L]$. Moreover, in view of Lemmas \ref{lem02} and \ref{lem1}, the sequence $n \mapsto \psi_{\epsilon,n}$ is increasing, and uniformly bounded on the closed annuli $\overline{B(0,L)}\setminus B(0,\rho)$, for every $\rho\in (1,L)$. By elliptic estimates \cite[\S 3.4]{gilbarg}, the first and second derivatives of $\psi_{\epsilon,n}$ are also uniformly bounded on these annuli. Thus, we can apply the theorem of Ascoli to $\psi_{\epsilon,n}$, and passing to the limit we obtain that $\psi_\epsilon:=\lim_{n\to\infty}\psi_{\epsilon,n}$ is a radial solution of \eqref{Dir2}. Setting, $\psi_\epsilon(x)=\tilde\psi_\epsilon(|x|)$, we have
\begin{equation}\label{polar}
\tilde \psi''_\epsilon(r)+\frac{n-1}{r}\tilde\psi'_\epsilon(r) =\epsilon^2W'(\tilde\psi_{\epsilon}(r)) \text{ on $(1,L]$},
\end{equation}  with $\lim_{r\to1}\tilde\psi_{\epsilon}(r)=\infty$, $\tilde\psi(L)=0$, and $\psi'_\epsilon<0$ on $(1,L]$. We also note that $\epsilon_1<\epsilon_2$ implies that $\psi_{\epsilon_1}\geq \psi_{\epsilon_2}$. Indeed, by Lemma \ref{lem02} we have $\psi_{\epsilon_1,n}\geq \psi_{\epsilon_2,n}$ for every $n$, since $\Delta \psi_{\epsilon_1,n}\leq \epsilon_2^2W'(\psi_{\epsilon_1,n})$.

It remains to show (ii) i.e. that $\lim_{\epsilon\to 0}\tilde\psi_\epsilon (r )=\infty$, holds for every $r \in (1,L)$ fixed. Assume by contradiction that $J:=\{\mu\in (1,L): \sup_\epsilon \tilde \psi_\epsilon(\mu)<\infty\}\neq\emptyset$. Our first claim is that $J$ is an open interval $(\mu_0,L)$, with $\mu_0\in [1,L)$. On the one hand it is clear that $\mu\in J$, implies that $[\mu,L)\subset J$. On the other hand, if $\mu\in J$, we obtain by elliptic estimates and the convexity of $\tilde \psi_\epsilon$, that $\tilde \psi'_\epsilon (L)$ is uniformly bounded, since $\tilde \psi_\epsilon$ is uniformly bounded on $[\mu,L]$. Next, an integration of \eqref{polar} over the interval $[a,b]\subset(1,L]$ gives
\begin{equation}\label{polar2}
\tilde \psi'_\epsilon(b)-\tilde \psi'_\epsilon(a)+(n-1)\frac{\tilde \psi_\epsilon(b)}{b}-(n-1)\frac{\tilde \psi_\epsilon(a)}{a}+(n-1)\int_a^b\frac{\tilde \psi_\epsilon(r)dr}{r^2} =\epsilon^2\int_a^b W'(\tilde\psi_{\epsilon}(r))dr.
\end{equation} 
Thus, setting $S:=\sup_\epsilon\tilde\psi_\epsilon(\mu)<\infty$, and defining $\rho_\epsilon\in (1,L)$ such that $\tilde \psi_\epsilon(\rho_\epsilon)=2S$, we deduce from \eqref{polar2} that the functions $|\tilde\psi'_\epsilon|$ are uniformly bounded on $[\rho_\epsilon,L]$ by a constant $M$. In particular we have $S\leq\int_{\rho_\epsilon}^\mu|\tilde\psi'_\epsilon|\leq M(\mu-\rho_\epsilon)$, and thus $\mu- \rho_\epsilon\geq\frac{S}{M}$. This proves our claim that $J$ is an open interval $(\mu_0,L)$. Since on the intervals $[\mu,L]$ with $\mu>\mu_0$, the functions $\tilde \psi_\epsilon$ are uniformly bounded up to the second derivatives, we can apply again the theorem of Ascoli, and obtain at the limit $\psi_0:=\lim_{\epsilon\to0}\psi_{\epsilon}$, where $\psi_0$ is a radial harmonic function defined on $\overline{B(0,L)}\setminus \overline{B(0,\mu_0)} $. Finally, we notice that $\lim_{|x|\to 1}\psi_0(x)=\infty$, since otherwise we would have $\mu_0\in J$. To conclude, we recall that a radial harmonic function is up to a multiplicative factor and an additive constant, the fundamental solution of Laplace's equation. Thus, the asymptotic condition $\lim_{|x|\to 1}\psi_0(x)=\infty$ cannot be realized.
\qed
\section{Proof of Theorem \ref{th1}}

\emph{Proof of (i). }
We consider the solution of the Dirichlet problem
\begin{equation}\label{Dir1a}
\Delta u_n =W'(u_n) \text{ in $\Omega_n:=\Omega\cap B(0,n)$, such that $u_n=n$ on $\partial \Omega_n$.}
\end{equation}
It is clear in view of Lemmas \ref{lem02} and \ref{lem1} that $u_n\geq 0$, $\forall n$, and that the sequence $u_n$ is uniformly bounded on compact subsets of $\Omega$. Moreover, by elliptic estimates \cite[\S 3.4]{gilbarg}, we obtain that the first and second derivatives of $u_n$ are as well  uniformly bounded on compact subsets. Thus, the theorem of Ascoli via a diagonal argument, implies that up to a subsequence, $u_n$ converges in $C^1_{\mathrm{loc}}(\Omega)$ to a solution of \eqref{pde} defined in $\Omega$. Next, since $\partial \Omega$ is Lipschitz \cite[\S 4.2]{evans-gariepy}, given $x_0\in \Omega$ and $\epsilon >0$, there exist $L>4$ independent of $\epsilon$, and an annulus $\omega_\epsilon:=B(y_\epsilon,L\epsilon)\setminus \overline{B(y_\epsilon,\epsilon)}$, such that $\overline{B(y_\epsilon,2\epsilon)}\subset \R^n\setminus \overline \Omega$, and $x_0\in B(y_\epsilon,L\epsilon/2)$. Finally, given $M>0$, we choose $\epsilon$ such that $\tilde\xi_\epsilon (L\epsilon/2)\geq M$ (cf. Lemma \ref{lem2} (ii)). Then, we see by Lemma \ref{lem02} that $n\geq \tilde \xi_\epsilon(2\epsilon)$ implies that $u_n\geq \xi_\epsilon(x-y_\epsilon)$ in $\omega_\epsilon \cap \Omega$, since $u_n=n\geq \xi_\epsilon(x-y_\epsilon)$ on $\overline\omega_\epsilon \cap \partial \Omega$, 
and $u_n\geq 0=\xi_\epsilon(x-y_\epsilon)$ on $\partial B(y_\epsilon,L\epsilon)\cap \overline\Omega$.
As a consequence, $u_n\geq \tilde \xi_\epsilon(L\epsilon/2)\geq M$ holds on $ B(y_\epsilon,L\epsilon/2)\cap \Omega$, for every $n\geq \tilde \xi_\epsilon(2\epsilon)$, and letting $n\to\infty$, we deduce that $u\geq M$ holds in $ B(y_\epsilon,L\epsilon/2)\cap \Omega$.
This proves that $\lim_{\Omega\ni x \to x_0} u(x)=+\infty$.
To complete the proof it remains to show that $u>0$. By construction it is obvious that $u \geq 0$. On the other hand, if $u(x_0)=0$ for some $x_0 \in \Omega$, then the maximum principle would imply that $u\equiv 0$, which is a contradiction.
\qed

\bigskip

\emph{Proof of (ii). }
Now we consider the solution of the Dirichlet problem
\begin{equation}\label{Dir1a}
\Delta u_n =W'(u_n) \text{ in $\Omega_n:=\Omega\cap B(0,n)$, such that $u_n=\phi_n$ on $\partial \Omega_n$,}\nonumber
	\end{equation}
with $\phi_n(x)=n\frac{d(x,A^-)-d(x,A^-)}{d(x,A^+)+d(x,A^-)}$. Clearly, $\phi_n$ is Lipschitz in $\Omega_n$, thus $\phi_n\in W^{1,2}(\Omega_n)$. In view of Proposition \ref{prop1} (a) and Lemma \ref{lem1} we have on the one hand 
\begin{equation}\label{lowb}
d(x,A^-\cup\partial B(0,n))>l(\lambda) \text{ for some $\lambda<0$} \Rightarrow u_n(x)\geq \lambda,
\end{equation}
and on the other hand
\begin{equation}\label{upb}
d(x,A^+\cup\partial B(0,n))>l(\lambda) \text{ for some $\lambda>0$} \Rightarrow u_n(x)\leq \lambda.
\end{equation}
Then, proceeding as in the proof of (i) we can see that $u_n$ converges in $C^1_{\mathrm{loc}}(\Omega)$ to a solution of \eqref{pde} defined in $\Omega$. To establish the boundary condition \eqref{bc3} in a neighbourhood of a point $x_0\in A^+$, we consider instead of $x\mapsto\xi_\epsilon(x-y_\epsilon)$ the barrier function $x\mapsto\xi_\epsilon(x-y_\epsilon)+\lambda$, where $\lambda<0$ is such that $d(x_0,A^-)>2l(\lambda)$. In view of \eqref{lowb}, the bound $u_n\geq \lambda$ holds on $B(x_0,l(\lambda))\cap \Omega$, uniformly in $n$. Finally, we notice that Lemma \ref{lem02} can be applied as previously to $u_n$ and $x\mapsto\xi_\epsilon(x-y_\epsilon)+\lambda$, since the latter function is a subsolution. In this way we establish that $\lim_{\Omega\ni x \to x_0} u(x)=+\infty$. The proof of \eqref{bc3} in a neighbourhood of a point $x_0\in A^-$ is similar.
\qed 

\bigskip

\emph{Proof of (iii). }
\begin{lemma}\label{proof3}
	Assuming \eqref{hhhh} and that the function $\Omega\ni x\mapsto d(x,\partial\Omega)$ is unbounded, then every solution $v$ of \eqref{pde} defined in $\Omega$ satisfies $\lim_{x\in\Omega,\ d(x,\partial \Omega)\to\infty} v(x)=0$.	
\end{lemma}
\emph{Proof. }
In view of Proposition \ref{prop1} (a) and Lemma \ref{lem1}, the condition $d(x,\partial\Omega)>l(\lambda)$ implies that $v(x)\leq \lambda$ (resp. $v(x)\geq \lambda$) for every $\lambda>0$ (resp. $\lambda<0$). Thus, the lemma follows by letting $\lambda\to 0$.\qed 

\bigskip

\emph{Proof of (iv) in the case of boundary condition \eqref{bc1}. }

Let $u$ and $v$ be two solutions of \eqref{pde} satisfying \eqref{bc1}. Our claim is that given $x_0\in \partial\Omega$, there exist $\epsilon>0$ and a constant $k>0$ such that $u+k\geq v$ holds for $x\in\Omega\cap B(x_0,\epsilon)$. Indeed, since $\partial \Omega$ is Lipschitz, there exists a ball $B(x_0,r)$ such that after a change of coordinates we have $B(x_0,r)\cap\Omega=\{x=(x_1,\ldots,x_n)\in B(x_0,r): x_n<f(x_1,\ldots,x_{n-1}) \}$, where $f:\R^{n-1}\to\R$ is a Lipschitz function \cite[\S 4.2]{evans-gariepy}. By taking $r$ smaller if necessary we may also assume that $u> 0$ on  $B(x_0,r)\cap\Omega$. Next we choose $\lambda>0$ such that $l(\lambda)=r/2$, and define $\omega_0:=B(x_0,r/2)\cap \Omega$. In view of Lemma \ref{lem1} and the convexity of $W'$ on $(0,\infty)$, it follows that $\omega_0\ni x\mapsto \chi(x)=u(x)+\tilde \alpha_\lambda (|x-x_0|)$ is a supersolution of \eqref{pde}. Moreover, we have that $\chi(x)\to\infty$, as $d(x,\partial \omega_0)\to 0$, with $x\in\omega_0$. Finally, Lemma \ref{lem02} applied in $\omega_0$ to $\chi$ and $x\mapsto v(x_1-\eta,x_2,\ldots,x_n)$ (with $\eta>0$ small), implies that $v(x_1-\eta,x_2,\ldots,x_n)\leq \chi(x)$, $\forall x\in \omega_0$. Letting $\eta\to 0$, we deduce that $v(x_1,x_2,\ldots,x_n)\leq \chi(x)$, $\forall x\in \omega_0$, thus our claim is established by taking $\epsilon=r/4$, and $k=\tilde \alpha_\lambda(r/4)$.

Clearly, when $\partial \Omega$ is compact, there exits a constant $k$ such that the inequality $v\leq u+k$ holds in a neighbourhood of $\partial\Omega$, and in addition $u(x)\to\infty$ (resp. $v(x)\to\infty$), as $d(x,\partial \Omega)\to 0$ (with $x\in\Omega$). 
As a consequence $u$ (resp. $v$) cannot take negative values in $\Omega$. Indeed, otherwise $u$ (resp. $v$) would attain in view of Lemma \ref{proof3} its negative minimum at a point $x_0\in \Omega$, and we would get $0\leq \Delta u(x_0)=W'(u(x_0))\in (0,\infty)$.

Next, given $\kappa>1$ and $\epsilon>0$, we consider the inequality $v\leq \kappa u+\epsilon$ holding when $d(x,\partial \Omega)\leq\eta$, with $\eta>0$ small. It follows from Lemma \ref{proof3} that we also have $v\leq \kappa u+\epsilon$ when $|x|\geq R$, with $R$ large enough. To conclude, we notice that $ \kappa u+\epsilon$ is a supersolution, since $u\geq 0$, and $W'$ is convex on $(0,\infty)$. Thus, Lemma \ref{lem02} implies that the inequality $v\leq \kappa u+\epsilon$ holds as well when $|x|< R$, and $d(x,\partial \Omega)>\eta$. This establishes that $v\leq \kappa u+\epsilon$ in $\Omega$. By letting $\kappa\to 1$, and $\epsilon\to 0$, we obtain $v\leq u$ in $\Omega$, and by interchanging $u$ with $v$ we get $u\equiv v$.
\qed 

\bigskip

\emph{Proof of (iv) in the case of boundary condition \eqref{bc3}. }

Now we consider $u$ and $v$ two solutions of \eqref{pde} satisfying \eqref{bc3}.
Since the subsets $A^\pm$ are compacts, we can show as previously that there exists a constant $k$ such that 
\begin{equation}\label{ascv}
\text{$|u-v|\leq k$ holds in a neighbourhood of $A^+$ (resp. $A^-$)}.
\end{equation}
In addition, we still have that
\begin{equation}\label{ascvv}
u(x),v(x)\to\infty \text{ (resp. $u(x),v(x)\to-\infty$), as $d(x,A^+)\to 0$ (resp. $d(x,A^-)\to 0$).}
\end{equation}

In the case of boundary condition \eqref{bc3}, the construction of the comparison function is much more involved. Let $\eta>0$. Our first claim is that
\begin{equation}\label{potW}
W'(t+\eta)\geq \frac{1}{2}W'(t), \ \forall t \in [-2\eta,0].
\end{equation}
This is obvious when $t\in[-\eta,0]$. On the other hand, if $t\in [-2\eta,\eta)$ we have $t+\eta\in [-\eta,0)$, and $\kappa:=\frac{t}{t+\eta}\geq 2$. By the concavity of $W'$ on $(-\infty,0)$, it follows that $W'(t)=W'(\kappa(t+\eta))\leq\kappa W'(t+\eta)\leq 2 W'(t+\eta)$, which establishes \eqref{potW}. 

Next, we define the constants:
\begin{itemize}
	\item $R$ such that $|x|\geq R\Rightarrow |u(x)| \text{ and } |v(x)|<\frac{\eta}{2}, \text{(cf. Lemma \ref{proof3})}$,
	\item $M=\sup\{ |\nabla u(x)|^2: \, |x|\leq R \text{ and } u(x)\in[-2\eta,0]\}$,
	\item $\mu \in \Big(0,\min\Big(\frac{1}{4M}, \frac{1}{4\int_{-\eta}^0(-W') } \Big)\Big)$,
	\item $\delta=\mu \int_{-\eta}^0(-W')\in \big(0,\frac{1}{4}\big)$,
\end{itemize}
and the continuous function
\begin{equation*}\label{fifi}
\phi(t)=\begin{cases}
0 &\text{ for } t\geq 0, \\
-\mu W'(t) &\text{ for } -\eta\leq t\leq 0, \\
\phi(-t-2\eta) &\text{ for } t\leq -\eta.
\end{cases}
\end{equation*}
Note that $\phi(t)\leq -\mu W'(t)$, for $t\leq 0$, since $W'$ is increasing on $\R$.
Let also $\Phi$ be the primitive of $\phi$ such that
\begin{equation*}\label{FiFi}
\Phi(t)=\int_{-\infty}^t\phi(s) ds +(1-\delta)=\begin{cases}
1+\delta &\text{ for } t\geq 0, \\
1 &\text{ for } t=-\eta, \\
1-\delta &\text{ for } t\leq -2\eta, \\
\end{cases}
\end{equation*}
and $f$ the primitive of $\Phi$ such that
\begin{equation*}\label{FiFi}
f(t)=\int_{0}^t\Phi(s) ds +2\eta=\begin{cases}
2\eta+(1+\delta)t &\text{ for } t\geq 0, \\
(1-\delta)(t+2\eta) &\text{ for } t\leq -2\eta.
\end{cases}
\end{equation*}
It is clear that $f$ is a smooth, convex, and strictly increasing function. One can easily check that
\begin{equation}\label{fprop}
f(t)>\eta+t, \ \forall t\in \R.
\end{equation}

We are going to show that $v\leq f(u)$ on $\Omega$. We first notice that the function $\Omega\ni x \mapsto \psi(x)=v(x)-f(u(x))$ satisfies $\psi(x)<0$ for $d(x,A^\pm)\leq \epsilon$, with $\epsilon>0$ small (cf. \eqref{ascv}, \eqref{ascvv}), and also for $|x|\geq R$ (cf. \eqref{fprop} and the definition of $R$).
To deduce that $\psi (x)\leq 0$ holds as well when $d(x,A^\pm)> \epsilon$, and $|x|< R$, it remains in view of Lemma \ref{lem02} to establish that $x \mapsto f(u(x))$ is a supersolution of \eqref{pde} in $\Omega \cap B(0,R)$. 

To see this, we distinguish the three following cases. If $u>0$, we have
\begin{align*}
\Delta (f(u))=(1+\delta)W'(u)\leq W'((1+\delta)u)\leq W'(2\eta+(1+\delta)u)=W'(f(u)),
\end{align*}
by the convexity of $W'$ on $(0,\infty)$, and the convexity of $W$. Otherwise, if $u<-2\eta$, we check as well that
\begin{align*}
\Delta (f(u))=(1-\delta)W'(u)\leq W'((1-\delta)u)\leq W'((1-\delta)(u+2\eta))=W'(f(u)),
\end{align*}
by the concavity of $W'$ on $(-\infty,0)$, and the convexity of $W$. Finally, if $-2\eta\leq u \leq 0$, we compute
\begin{align*}
\Delta (f(u))&=f'(u)W'(u)+f''(u)|\nabla u|^2,
\end{align*}
and we have
\begin{itemize}
	\item $f'(u)W'(u)\leq \frac{3}{4} W'(u)$, since $f'\geq 1-\delta\geq \frac{3}{4}$.
	\item $f''(u)|\nabla u|^2=\phi(u)|\nabla u|^2 \leq -\mu M W'(u)\leq -\frac{1}{4}W'(u)$. 
\end{itemize}
Thus, we still get $\Delta (f(u))\leq \frac{1}{2}W'(u)\leq W'(u+\eta)\leq W'(f(u))$, by \eqref{potW}, \eqref{fprop}, and the convexity of $W$. 

This proves that $v\leq f(u)$ holds in $\Omega$, for every choice of $\eta>0$, and $\mu \in \Big(0,\min\Big(\frac{1}{4M}, \frac{1}{4\int_{-\eta}^0(-W') } \Big)\Big)$. By letting $\eta\to 0$, and $\mu \to 0$, it follows that $v\leq u$ holds in $\Omega$, and by interchanging $v$ with $u$, we conclude that $u\equiv v$. \qed

\section*{Acknowledgments}
The author was partially supported by the National Science Centre, Poland (Grant No. 2017/26/E/ST1/00817).
He would also like to thank Anestis Fotiadis for several fruitful discussions that motivated the study of this problem.

\bibliographystyle{plain}

\end{document}